\documentclass[11pt]{article} 
\usepackage{amsmath,amssymb,amsfonts,amsthm,graphicx}
\newtheorem{theorem}{Theorem}

\newtheorem{lem}{Lemma}
\newtheorem{defin}{Definition}

\title{Optimal Control for a Nonlinear Stochastic Parabolic Model of Language Competition}
\author{ {$^{a}$ Sakine Esmaili\footnote{{\em E-mail addresses:}~sakine.esmaili@modares.ac.ir
(Sakine Esmaili), eslahchi@modares.ac.ir (M. R. Eslahchi).} \hspace{.5cm} $^{a}$ M. R. Eslahchi}
\vspace{.5cm}$$\\
\small{\em $$\em $^{a}$ Department of Applied Mathematics, 
Faculty of Mathematical Sciences,}\vspace{-1mm}\\
\small{\em Tarbiat Modares University, P.O. Box 14115-134,
Tehran, Iran}}
\makeatletter 
\def\smallunderbrace#1{\mathop{\vtop{\m@th\ialign{##\crcr 
$\hfil\displaystyle{#1}\hfil$\crcr 
\noalign{\kern3\p@\nointerlineskip}%
\tiny\upbracefill\crcr\noalign{\kern3\p@}}}}\limits} 

\makeatother 
   
\begin{document}

\maketitle 
\begin{abstract}
In this investigation, an optimal control problem for a stochastic mathematical model of language competition is studied. We have considered  the stochastic model of language competition by adding the stochastic terms to the deterministic model to take into account  the random perturbations and uncertainties caused by the environment to have more reliable model.  The model has formulated the population densities of the speakers of two languages, which are competing against each other  to be saved from destruction, attract more speakers and so on, using two nonlinear stochastic parabolic equations.  Four factors including the status of the languages and the growth rates of the populations are considered as the control variables (which can be controlled by the speakers of the populations or policy makers who make decisions for the populations for particular purposes)  to control the evolution of the population densities. Then,  the optimal control problem for the stochastic model of language competition is studied. Employing the tangent-normal cone techniques, the Ekeland variational principle and other theorems proved throughout the paper, we have shown  there exists unique
stochastic optimal control. We have also presented the exact form of the optimal control in terms of stochastic adjoint states. 

\end{abstract}
\textbf{Keywords:} {Optimal control,  Nonlinear stochastic parabolic equation,  Existence and uniqueness of stochastic optimal control, Language competition.}\\

\noindent\textbf{MSC:} {49K20, 49J20, 35R60. }
\section{Introduction}
There are different languages spoken by different populations,  which are affected by a variety of factors consisting of growth rates of the populations, status of the languages and movement of the speakers,  and may be in danger of extinction. It seems that there is  a competition between the languages to  attract more speakers, increase their population sizes and
be saved from extinction, especially when different languages are spoken in the same area.  For this reason, it is of great importance to study  the evolution  of  endangered languages to be able to save them from extinction.  Hence, the mathematical modelling of language competition is applied to study the evolution of speaker populations of languages   [1--8], which leads us to predict the influence  of different factors on the future of  them and detect the efficient factors  in saving at-risk languages. For instance, since India, with a 125 million English speakers,
has the second largest number of English speakers in the world, the model of interaction
between Monolingual Hindi speakers, Hindi/English Bilinguals and Hinglish speakers is studied in \cite{parshad}. 
  In another paper  \cite{iri}, the authors  have investigated a society with two official languages consisting of A,  spoken by every individual in the society, and
B,  spoken by a bilingual minority. In \cite{Bakalis}, time
evolution of the density of speakers of an at-risk language, Aromanian, which is spoken by a
bilingual community in North-Western Greece, is studied.  In investigation \cite{patrica},  the  language competition with game theoretic approach  is analyzed when one  language is known
by all individuals and  the other one is  spoken by a minority. The author of \cite{wiese},   has also used game theory  to contribute to  the economics of language and  offered some explanations of how  past languages like
Latin or Sanskrit can develop into a standard for literary production. In another investigation   \cite{lan1}, the authors  studied a mathematical model of language competition consisting of two coupled nonlinear parabolic
equations. Each parabolic equation presents the population density of the speakers of one language,
which is competing against another one to be saved from destruction. 
In this model, the effects of some factors including growth rates of the populations and status of the languages on the evolution of the population densities are modelled. Since the language with lower status or the population with lower growth rate may be destroyed, one may think that by changing these factors, we can  save the languages. But in reality, every change can have negative effect on something else and  can be  so costly. Therefore, considering optimal control for  the model is of great importance and helps us to adopt accurate policies for saving at-risk languages.   Some other examples of competition between populations can be seen in \cite{phys1,mar1,bac}.\\
Owing to the fact that all systems are affected by many uncertainties  caused by the environment,  considering the models in which the stochastic terms are included, enables us to deal with more trustable models.  In this paper, we have studied a stochastic optimal bounded control strategy for a nonlinear  stochastic parabolic model of language competition with a quadratic cost function. We have considered the mathematical model of language competition studied in \cite{lan1} but we have added stochastic terms to the model to consider uncertainties and random perturbations and study more reliable model.  In order to control  the population densities  of the speakers of two languages,  our control strategy is  considering some factors in the nonlinear stochastic parabolic model  including the status of the languages and the growth rates of the populations  to play role as the control variables, which have direct effects on the languages to attract  more speakers.  Since changing the control variables unthinkingly in order to save  an at-risk language can be very costly and endanger  another one, we have considered an optimal control problem for the model to minimize a  quadratic  cost function in which the cost of changing the values of status and growth rates are considered and our goals for population sizes are included. Designing  this cost function helps us to make accurate and efficient decisions. For deriving the bounded optimal control variables, we have used some stochastic auxiliary equations (so-called  adjoint equations). Employing the stochastic adjoint equations, some variational techniques and the tangent-normal cone techniques, the necessary conditions for the stochastic optimal control  variables (which are  stochastic because they are obtained in terms of the stochastic adjoint equations) are presented. After that, applying the Ekeland variational principle and adjoint equations, we have made a sequence converging  to the exact optimal bounded control variables. Using the constructed sequence, the existence and uniqueness of the stochastic optimal control  are proved. Since the  control variables are  assumed to be bounded,  the obtained  optimal control variables are  continuous piecewise-defined functions of the stochastic adjoint states. \\
The organization of paper is as follows. In Section 2, the model and the optimal control problem are presented. The adjoint equations are introduced in Subsection 2.1. The existence and uniqueness of optimal control are studied in Subsection 2.2. Concluding remarks  are offered in Section 3. Some lemmas
and mathematical concepts, which are used throughout the paper are included in Appendix. 
\section{Model and  Optimal Control Problem}\label{OCP1111}
\noindent In this section, first we introduce the stochastic mathematical model of language competition  including two nonlinear stochastic parabolic equations describing the evolution of the population densities of the speakers of two languages. The stochastic model  is obtained by adding stochastic terms to the model studied in \cite{lan1} to have more reliable model by considering the random perturbations and uncertainties.
\begin{equation}\label{YAS1111}
{d{f_1}}=\nabla.(D_{1}\nabla{  {f_1}})dt-\nabla.(F_1f_1)dt+R_1(f_1,f_2,\beta_1)dt+R(f_1,f_2,s_1,s_2)dt+f_1\sum_{i=0}^nh_{i,1}(x,y)dB_i,
\end{equation}
\[~(x,y,t)\in\Omega\times (0,\infty),\]
\begin{equation}\label{313}
{f_1}\left({x,y,}t \right){\rm =}{f_1^b}\left({x,y,}t \right), ~(x,y)\in\partial\Omega,~ {t >}0,~{f_1}\left(x,y ,0\right){\rm =}{f^0_1}\left(x,y \right),~(x,y)\in\Omega,
\end{equation} 
\begin{equation}\label{YIHkhco}
{d {f_2}}=\nabla.(D_{2}\nabla{  {f_2}})dt-\nabla.(F_2f_2)dt+R_2(f_1,f_2,\beta_2)dt-R(f_1,f_2,s_1,s_2)dt+f_2\sum_{i=0}^nh_{i,2}(x,y)dB_i,
\end{equation}
\[~(x,y,t)\in\Omega\times(0,\infty),\]
\begin{equation}\label{YAS1*ekhco}
{f_2}\left({x,y,}t \right){\rm =}{f_2^b}\left({x,y,}t \right), ~(x,y)\in\partial\Omega,~ {t >}0,~{f_2}\left(x,y ,0\right){\rm =}{f^0_2}\left(x,y \right),~(x,y)\in\Omega,
\end{equation}
\[R_1(f_1,f_2,\beta_1)=\beta_1f_1(1-\dfrac{f_1+f_2}{N}),~~R_2(f_1,f_2,\beta_2)=\beta_2f_2(1-\dfrac{f_1+f_2}{N}),\]
\[R(f_1,f_2,s_1,s_2)=k(s_1f_1^\alpha f_2-s_2f_1f_2^\alpha),\]
where $B(t)=(B_1(t),\cdots,B_n(t))$ is an $n$-dimensional  Brownian motion on $(\mathbf{\Omega},\mathfrak{F},\mathbb{P})$ with associated filtration $\{\mathfrak{F}(t)\}_{t>0}$, ${f_i}\left({x,y,}t \right)$ describes the population density of speakers of language $i$, $\alpha$ is a positive constant,  $s_i(t)$ represents the status of
language $i$ at time $t$ and the logistic terms with  carrying capacity $N$ and Malthus rates $\beta_1$ and $\beta_2$  take into account the
growth of the populations. The population movement is described by the advection terms with  the external force fields $F_1(x, y)=(F_1^x(x,y),F_1^y(x,y))$ and $F_2(x, y)=(F_2^x(x,y),F_2^y(x,y))$
and by the diffusion terms with the diffusion coefficients $D_1$ and $D_2$.\\

\noindent In this model, we have assumed that the following properties are satisfied: \\
\noindent\textbf{A.} $\Omega\subset \mathbb{R}^2$ is an open, bounded set   with boundary $\partial\Omega\in C^2$.\\
\noindent\textbf{B.} There exist $C^2-$smooth functions  $u_1(x,y,t)$ and $u_2(x,y,t)$  such that $u_i(x,y,t)={f^b_i}\left(x,y ,t\right)$  $(i=1,2)$  on $\partial \Omega$ and $u_i(x,y,0)={f^0_i}\left(x,y\right)$ $(i=1,2)$ on $\Omega$.\\
\noindent\textbf{C.} $F_1$ and $F_2$  are $C^1-$smooth functions, $h_{i,j}$ ($i=0,\cdots,n,~j=1,2$) are  $C^2-$smooth functions and  $h_{i,j}(x,y)=0$ on $\partial \Omega$.

\noindent In the following theorem the existence and uniqueness of the solution of \eqref{YAS1111}--\eqref{YAS1*ekhco} are presented. 
\begin{theorem}\label{mainthm1}
 Let assumptions \textbf{A}--\textbf{C}   be satisfied, then the problem \eqref{YAS1111}--\eqref{YAS1*ekhco} for almost every $\omega\in\mathbf{\Omega}$ has a unique solution. For every $T>0$, $f_1$ and $f_2$ are $\mathfrak{F}(t)$-adapted and for almost every $\omega\in \mathbf{\Omega}$,  we also have
 \[{f_1},{f_2}\in C\left(\Omega\times [0,T]\right).\]
\end{theorem}
\noindent{\bf{Proof}} Using Lemma \ref{Lemma 1ito} (see Appendix) and fixed point theorems, we can prove this theorem.
\qed 
 
\noindent Now, we  consider  the optimal control for the presented  model of language competition.   In this language competition optimal control problem  (LCOCP), we purpose to control the population densities of the speakers by considering some time dependant control variables $\beta_1$, $\beta_2$, $s_1$, and $s_2$, which control the growth of the populations and the status of the languages. Then, we minimize
\[J(\beta_1,\beta_2,s_1,s_2):=\int_\Omega\int_0^T((\beta_1-b_1)^2+(\beta_2-b_2)^2+\lambda_1(f_1-r_1^*)^2+\lambda_2(f_2-r_2^*)^2)dt dxdy+\]
\begin{equation}\label{focpcostf}
\int_\Omega\int_0^T(s_1-r_3^*)^2+(s_2-r_4^*)^2dtdxdy,
\end{equation}
\[~~~~(\beta_1,\beta_2,s_1,s_2)\in \mathcal{C}_{ad}\times\mathcal{W}_{ad}\times\mathcal{C}^1_{ad}\times\mathcal{W}^1_{ad},\]
where $r^*_1,r^*_2,r_3^*,r_4^*,b_1,b_2\in C^{\alpha,\alpha/2}(\overline{\Omega}\times [0, T]) $  for some $0<\alpha<1$, also 
\begin{equation}\label{admis}
\mathcal{C}_{ad}:=\Big\{ v\in C(\overline{\Omega}\times [0, T])~~a.s.: ~v~\textup{is}~\mathfrak{F}(t)-\textup{adapted},~ l^{*}_1(x,y,t)\leq v(x,y,t)\leq l^{**}_1(x,y,t)~~a.s.\Big\},
\end{equation}
\begin{equation}\label{admiss}
 \mathcal{W}_{ad}:=\Big\{ v\in C(\overline{\Omega}\times [0, T])~~a.s.:  ~v~\textup{is}~\mathfrak{F}(t)-\textup{adapted},~ l^{*}_2(x,y,t) \leq v(x,y,t)\leq l^{**}_2(x,y,t)~~a.s.\Big\},
\end{equation}
\begin{equation}
\mathcal{C}^1_{ad}:=\Big\{ v\in C(\overline{\Omega}\times [0, T])~~a.s.: ~v~\textup{is}~\mathfrak{F}(t)-\textup{adapted},~ l^{*}_3(x,y,t)\leq v(x,y,t)\leq l^{**}_3(x,y,t)~~a.s.\Big\},
\end{equation}
\begin{equation}\label{admis11}
\mathcal{W}^1_{ad}:=\Big\{ v\in C(\overline{\Omega}\times [0, T])~~a.s.: ~v~\textup{is}~\mathfrak{F}(t)-\textup{adapted},~ l^{*}_4(x,y,t)\leq v(x,y,t)\leq l^{**}_4(x,y,t)~~a.s.\Big\},
\end{equation}
\[l_i^{*},l_i^{**}\in C^{\alpha,\alpha/2}(\overline{\Omega}\times [0, T]),~i=1,2,3,4,\]
such that the equations \eqref{YAS1111}--\eqref{YAS1*ekhco} are satisfied.

\subsection{ٍAdjoint Equations}\label{EUOP7777}
Before presenting our main results,  we introduce the following stochastic adjoint equations, which are instrumental in presenting the necessary conditions and proving the existence and  uniqueness of stochastic optimal control variables. The optimal control variables are stochastic because they will obtained in terms of the stochastic adjoint states.
\[d z_{f_1}=-\nabla.(D_{1}\nabla{  z_{f_1}})dt+\]
\begin{equation}\label{YAS*1e1}
H^*_1(\smallunderbrace{z_{f_1},z_{f_2},f_1,f_2,\beta_1,\beta_2,s_1,s_2}_{Z^*})dt-z_{f_{1}}\sum_{i=0}^nh_{i,1}(x,y)dB_i,
\end{equation}
\[~(x,y,t)\in\Omega\times (0,T),\]
\begin{equation}
z_{f_1}\left({x,y,}t \right){\rm =}0, ~(x,y)\in\partial\Omega,~ 0<t< T,~z_{f_1}\left(x,y ,T\right){\rm =}0,~(x,y)\in\Omega, 
\end{equation} 

\begin{equation}\label{YIH}
d z_{f_2}=-\nabla.(D_{2}\nabla{  z_{f_2}})dt+H^*_2\left({Z^*}\right)dt-z_{f_{2}}\sum_{i=0}^nh_{i,2}(x,y)dB_i,~(x,y,t)\in\Omega\times (0,T),
\end{equation}
\begin{equation}\label{YAS1*e}
z_{f_2}\left({x,y,}t \right){\rm =}0, ~(x,y)\in\partial\Omega,~ 0< t< T,~z_{f_2}\left(x,y ,T\right){\rm =}0,~(x,y)\in\Omega,
\end{equation}
\begin{eqnarray}\label{pqd*e}
&d z_{g_1}=
G_1^*\left(Z^*\right)dt,~d z_{g_2}= 
G_2^*\left(Z^*\right)dt,~(x,y)\in\Omega,~ 0< t< T,&\\  
&d z_{g_3}{\rm =}
G_3^*\left(Z^*\right)dt,~d z_{g_4}{\rm =}
G_4^*\left(Z^*\right)dt,~(x,y)\in\Omega,~0< t< T,&
\end{eqnarray}
\begin{eqnarray}\label{HM345**ee}
&z_{g_i}\left(x,y,T\right)=0, ~(x,y)\in\Omega,~i=1,2,3,4,&
\end{eqnarray}
where 
\[H^*_1\left(Z^*\right)=-\dfrac{\partial \Big(R_1(f_1,f_2,\beta_1)+ R(f_1,f_2,s_1,s_2)\Big)}{\partial f_1}z_{f_1}-\]
\[\dfrac{\partial \Big(R_2(f_1,f_2,\beta_2)-R(f_1,f_2,s_1,s_2)\Big)}{\partial f_1}z_{f_2}-\nabla(z_{f_1}).F_1+\lambda_1(f_1-r_1^*)+z_{f_{1}}\sum_{i=0}^n(h_{i,1}(x,y))^2,\]
\[H^*_2\left(Z^*\right)=-\dfrac{\partial \Big(R_1(f_1,f_2,\beta_1)+ R(f_1,f_2,s_1,s_2)\Big)}{\partial f_2}z_{f_1}-\]
\[\dfrac{\partial \Big(R_2(f_1,f_2,\beta_2)-R(f_1,f_2,s_1,s_2)\Big)}{\partial f_2}z_{f_2}-\nabla(z_{f_2}).F_2+\lambda_2(f_2-r_2^*)+z_{f_{2}}\sum_{i=0}^n(h_{i,2}(x,y))^2,\]
\begin{equation}\label{1122}
G^*_1\left(Z^*\right)=-\dfrac{\partial R_1(f_1,f_2,\beta_1)}{\partial \beta_1}z_{f_1},
\end{equation}
\begin{equation}
G^*_2\left(Z^*\right)=-\dfrac{\partial R_2(f_1,f_2,\beta_2)}{\partial \beta_2}z_{f_2},
\end{equation}
\begin{equation}
G^*_3\left(Z^*\right)=-\dfrac{\partial R(f_1,f_2,s_1,s_2)}{\partial s_1}z_{f_1}+\dfrac{\partial R(f_1,f_2,s_1,s_2)}{\partial s_1}z_{f_2},
\end{equation}
\begin{equation}\label{1113}
G^*_4\left(Z^*\right)=-\dfrac{\partial R(f_1,f_2,s_1,s_2)}{\partial s_2}z_{f_1}+\dfrac{\partial R(f_1,f_2,s_1,s_2)}{\partial s_2}z_{f_2}.
\end{equation}
\noindent In the following theorem, the existence and uniqueness of the solution of   the system of  stochastic adjoint equations are presented.

\begin{theorem}\label{AMHZ}
Let assumptions \textbf{A}--\textbf{C} be satisfied, then for almost every $\omega\in\mathbf{\Omega}$, the adjoint system  \eqref{YAS*1e1}--\eqref{HM345**ee} has a unique  solution $(z_{f_{1}},z_{f_{2}},z_{g_{1}},z_{g_{2}},z_{g_{3}},z_{g_{4}})$ and for every $T>0$,  we have   $z_{f_{1}}$, $z_{f_{2}}$, $z_{g_{1}}$, $z_{g_{2}}$, $z_{g_{3}}$ and $z_{g_{4}}$ are $\mathfrak{F}(t)$-adapted and
 \[z_{f_{1}},z_{f_{2}},z_{g_{1}},z_{g_{2}},z_{g_{3}},z_{g_{4}}\in C\left(\Omega\times [0,T]\right). ~~a.s.\]
\end{theorem}
\noindent{\bf{ Proof}} Employing Lemma \ref{Lemma 1ito} (see Appendix) and Theorem \ref{mainthm1},  we  obtain the  result.
\qed
\subsection{Existence and Uniqueness of Optimal Control}\label{EUOP1111}
 In this subsection, we present the necessary conditions for optimal control, which give us the explicit forms of stochastic optimal control variables in terms of the stochastic adjoint states, also we prove the existence and  uniqueness of optimal control variables. \\
   In the next theorem, we introduce the necessary conditions for optimal control variables of  LCOCP.
\begin{theorem}\label{ncon} \textup{\textbf{(Necessary Conditions)}}
Let $(\beta_1^*,\beta^*_2,s_1^*, s_2^*)$ be the optimal control of LCOCP, then for almost every $\omega\in \mathbf{\Omega}$, we have
\begin{equation}\label{YA1**e5e}
\beta_1^*=\mathcal{F}\Big(b_1-G_1^*({Z^*}),l_1^{*},l_1^{**}\Big),~
\beta_2^*=\mathcal{F}\Big(b_2-G_2^*\left(Z^*\right),l_2^{*},l_2^{**}\Big),
\end{equation}
\begin{equation}\label{YA1**e}
~s_1^*=\mathcal{F}\Big(r_3^*-G_3^*\left(Z^*\right),l_3^*,l_3^{**}\Big),~s_2^*=\mathcal{F}\Big(r_4^*-G_4^*\left(Z^*\right),l_4^*,l_4^{**}\Big),
\end{equation}
where  $(z_{f_1},z_{f_2},f_1,f_2)$ is the exact solution of \eqref{YAS1111}--\eqref{YAS1*ekhco} and \eqref{YAS*1e1}--\eqref{YAS1*e}  corresponding to $(\beta^*_1,\beta^*_2,s^*_1, s^*_2)$,   $\beta^*_1$, $\beta^*_2$, $s^*_1$ and $ s^*_2$  are $\mathfrak{F}(t)$-{adapted},~ $l_i^{*}$ and $l_i^{**}~(i=1,2,3,4)$ are defined in \eqref{admis}-\eqref{admis11}, and 
\begin{equation}\label{YA2**}
\mathcal{F}(x,u,v)=\begin{cases}
u,~~x<u,\\
x,~~u\leq x\leq v,\\
v,~~v< x .\\
\end{cases}
\end{equation}
\end{theorem}
\noindent{ \textbf{Proof} Let $\left(f_1^*,f_2^*\right)$ and $\left(f_1^{\epsilon},f_2^{\epsilon}\right)$ be the solutions of \eqref{YAS1111}--\eqref{YAS1*ekhco} corresponding to the controls  \[(\beta_1^*,\beta^*_2,s_1^*, s_2^*) \in  \mathcal{C}_{ad}\times\mathcal{W}_{ad}\times\mathcal{C}^1_{ad}\times\mathcal{W}^1_{ad},\]
 and 
 \[(\beta_1^\epsilon,\beta^\epsilon_2,s_1^\epsilon, s_2^\epsilon) =(\beta_1^*+\epsilon \beta_1^0,\beta_2^*+\epsilon \beta_2^0,s_1^*+\epsilon s_1^0,s_2^*+\epsilon s_2^0) \in  \mathcal{C}_{ad}\times\mathcal{W}_{ad}\times\mathcal{C}^1_{ad}\times\mathcal{W}^1_{ad},\]
respectively, where $\mathcal{C}_{ad}$,  $\mathcal{W}_{ad}$, $\mathcal{C}^1_{ad}$ and  $\mathcal{W}^1_{ad}$ are defined in \eqref{admis}--\eqref{admis11} and $\epsilon$  is a positive constant. Also, assume that  
\[J(\beta_1^*,\beta^*_2,s_1^*, s_2^*)=\min \{J(\beta_1,\beta_2,s_1, s_2):(\beta_1,\beta_2,s_1,s_2)\in \mathcal{C}_{ad}\times\mathcal{W}_{ad}\times\mathcal{C}^1_{ad}\times\mathcal{W}^1_{ad}\},\]
\[f_{1,1}^{\epsilon}=\dfrac{f_1^{\epsilon}-f_1^*}{\epsilon},~f_{2,1}^{\epsilon}=\dfrac{f_2^{\epsilon}-f_2^*}{\epsilon}.\]
Therefore, one can derive that
\[\dfrac{J(\beta_1^\epsilon,\beta_2^\epsilon,s_1^{\epsilon},s_2^{\epsilon})-J(\beta_1^*,\beta_2^*,s_1^*,s_2^*)}{\epsilon}=\]
\[\int_\Omega\int_0^T((\beta^*_1+\beta_1^\epsilon-2b_1)\beta_{1}^0+(\beta^*_2+\beta_2^\epsilon-2b_2)\beta_{2}^0+\lambda_1(f^*_1+f_1^\epsilon-2r_1^*)f_{1,1}^\epsilon+\lambda_2(f_2^*+f_2^\epsilon-2r_2^*)f_{2,1}^\epsilon)dtdxdy+\]
\begin{equation}\label{AA*}
\int_\Omega\int_0^T(s_1^*+s^\epsilon_1-2r_3^*)s_1^0+(s_2^*+s^\epsilon_2-2r_4^*)s_2^0 dtdxdy\geq 0,~~a.s.
\end{equation}
and
\begin{equation}\label{YAS*1}
d{f_{1,1}^\epsilon}=\nabla.(D_{1}\nabla{  {f_{1,1}^\epsilon}})dt+E^\epsilon_1(\smallunderbrace{f_1^*,f_2^*,\beta_1^*,\beta_2^*,s_1^*,s_2^*,f_1^\epsilon,f_2^\epsilon,\beta_1^\epsilon,\beta_2^\epsilon,s_1^\epsilon,s_2^\epsilon}_{Z^{*,\epsilon}})dt+f^\epsilon_{1,1}\sum_{i=0}^nh_{i,1}dB_i,
\end{equation}
\[~(x,y,t)\in\Omega\times (0,\infty),\]
\begin{equation}
{f_{1,1}^\epsilon}\left({x,y,}t \right){\rm =}0, ~(x,y)\in\partial\Omega,~ {t >}0,~f_{1,1}^\epsilon\left(x,y ,0\right){\rm =}0,~(x,y)\in\Omega,
\end{equation} 
\begin{equation}\label{YIH}
d {f_{2,1}^\epsilon}=\nabla.(D_{2}\nabla{  {f_{2,1}^\epsilon}})dt+E^\epsilon_2({Z^{*,\epsilon}})dt+f^\epsilon_{2,1}\sum_{i=0}^nh_{i,2}dB_i,~(x,y,t)\in\Omega\times (0,\infty),
\end{equation}
\begin{equation}\label{YAS1*e11}
{f_{2,1}^\epsilon}\left({x,y,}t \right){\rm =}0, ~(x,y)\in\partial\Omega,~ {t >}0,~{f_{2,1}^\epsilon}\left(x,y ,0\right){\rm =}0,~(x,y)\in\Omega,
\end{equation}
where
\[E^\epsilon_1({Z^{*,\epsilon}})=-\nabla.(F_1 f_{1,1}^\epsilon)+\dfrac{R_1(f^\epsilon_1,f^\epsilon_2,\beta_1^\epsilon)+R(f^\epsilon_1,f^\epsilon_2,s^\epsilon_1,s^\epsilon_2)-\Big(R_1(f^*_1,f^*_2,\beta_1^*)+R(f^*_1,f^*_2,s^*_1,s^*_2)\Big)}{\epsilon},\]
\[E_2^\epsilon({Z^{*,\epsilon}})=-\nabla.(F_2 f_{2,1}^\epsilon)+\dfrac{R_2(f^\epsilon_1,f^\epsilon_2,\beta_2^\epsilon)-R(f^\epsilon_1,f^\epsilon_2,s^\epsilon_1,s^\epsilon_2)-\Big(R_2(f^*_1,f^*_2,\beta_2^*)-R(f^*_1,f^*_2,s^*_1,s^*_2)\Big)}{\epsilon}.\]
Therefore using \eqref{YAS*1e1}--\eqref{HM345**ee} and   Theorem \ref{??} (see Appendix), one can deduce that
\[\int_0^T\int_\Omega d(f^\epsilon_{1,1} z_{f_1})+d(f^\epsilon_{2,1} z_{f_2})+
\int_0^T\int_{\Omega}\beta_1^0d z_{g_1}+\beta_2^0d z_{g_2}+s_1^0dz_{g_3}+
s_2^0 d z_{g_4}dxdy=\]
\[\int_0^T\int_\Omega \Big(f^\epsilon_{1,1}d z_{f_1}+z_{f_1}df^\epsilon_{1,1}+dz_{f_1}df^\epsilon_{1,1}+f^\epsilon_{2,1}d z_{f_2}+z_{f_2}df^\epsilon_{2,1}+dz_{f_2}df^\epsilon_{2,1}\Big)+\]
\[\int_0^T\int_{\Omega}\Big(\beta_1^0d z_{g_1}+\beta_2^0dz_{g_2}+s_1^0d z_{g_3}+
s_2^0d z_{g_4}\Big) dxdy =\]
\[\int_0^T\int_\Omega\Big(\lambda_1(f_1^*-r_1^*)f^\epsilon_{1,1}+\lambda_2(f_2^*-r_2^*)f^\epsilon_{2,1}\Big)dxdy dt+\varpi(\epsilon, T),\]
where $(z_{f_1},z_{f_2},z_{g_1},z_{g_2},z_{g_3},z_{g_4})$ is the exact solution of the adjoint system \eqref{YAS*1e1}--\eqref{HM345**ee}  corresponding to $(\beta^*_1,\beta^*_2,s^*_1, s^*_2)$ and $\displaystyle{\lim_{\epsilon\rightarrow 0}}\varpi(\epsilon, T)=0$, thus we have
\begin{equation}\label{AA*e}
\int_0^T\int_\Omega\Big(\lambda_1(f_1^*-r_1^*)f^\epsilon_{1,1}+\lambda_2(f_2^*-r_2^*)f^\epsilon_{2,1}-\beta_1^0\dfrac{\partial z_{g_1}}{\partial t}-\beta_2^0\dfrac{\partial z_{g_2}}{\partial t}-s_1^0\dfrac{\partial z_{g_3}}{\partial t}-s_2^0\dfrac{\partial z_{g_4}}{\partial t}\Big)dxdy d t=-\varpi(\epsilon, T).
\end{equation}
Hence, using \eqref{AA*} and  \eqref{AA*e}, we arrive at
\[\int_\Omega\int_0^T(\beta^*_1-b_1+\dfrac{\partial z_{g_1}}{\partial t})\beta_{1}^0+(\beta^*_2-b_2+\dfrac{\partial z_{g_2}}{\partial t})\beta_{2}^0dtdxdy+\]
\begin{equation}\label{AA*ee}
\int_{\Omega}\int_0^T (s_1^*-r_3^*+\dfrac{\partial z_{g_3}}{\partial t})s_1^0+(s_2^*-r_4^*+\dfrac{\partial z_{g_4}}{\partial t})s_2^0dtdxdy\geq 0,~~a.s.
\end{equation}
Therefore using tangent-normal cone techniques \cite{25.} and \eqref{admis}--\eqref{admis11}, $(\beta_1^*,\beta^*_2,s_1^*, s_2^*)$ is as follows
\begin{equation}\label{habb}
\beta_1^*=\mathcal{F}\Big(b_1-\dfrac{\partial z_{g_1}}{\partial t},l_1^{*},l_1^{**}\Big),~\beta_2^*=\mathcal{F}\Big(b_2-\dfrac{\partial z_{g_2}}{\partial t},l_2^{*},l_2^{**}\Big),~~a.s.
\end{equation}
\begin{equation}\label{habb1}
~s_1^*=\mathcal{F}\Big(r_3^*-\dfrac{\partial z_{g_3}}{\partial t},l_3^*,l_3^{**}\Big),~s_2^*=\mathcal{F}\Big(r_4^*-\dfrac{\partial z_{g_4}}{\partial t},l_4^*,l_4^{**}\Big).~~a.s.
\end{equation}
Also from \eqref{habb}--\eqref{habb1}, we can conclude that   $\beta^*_1$, $\beta^*_2$, $s^*_1$ and $ s^*_2$  are $\mathfrak{F}(t)$-{adapted}.
\qed
\\
\noindent The following theorem gives the existence and  uniqueness of the optimal control variables for  LCOCP.
\begin{theorem}\label{YKHFKH} \textup{\textbf{(Existence and Uniqueness)}}
 LCOCP, for almost every $\omega\in \mathbf{\Omega}$, has a unique optimal control, which satisfies \eqref{YA1**e5e}--\eqref{YA1**e}.
\end{theorem}
\noindent{\textbf{Proof } Using Theorem \ref{??} (see Appendix),  we conclude that  for each 
\[(\beta_1,\beta_2,s_1, s_2)\in\mathcal{C}^*_{ad}\times\mathcal{W}^*_{ad}\times\mathcal{C}^{1*}_{ad}\times\mathcal{W}^{1*}_{ad}\setminus C_4(\smallunderbrace{\Omega\times [0, T]}_{\overline{Q_T}}),\]
where
\[C_4({\overline{Q_T}})=C(\overline{Q_T})\times C(\overline{Q_T})\times C(\overline{Q_T})\times C(\overline{Q_T}),\]
\begin{equation*}
\mathcal{C}^*_{ad}:=\Big\{ v\in L^\infty(\overline{Q_T})~~a.s.: l^{*}_1(x,y,t) \leq v(x,y,t)\leq l^{**}_1(x,y,t),~~a.e.~~a.s. \Big\},
\end{equation*}
\begin{equation*}
 \mathcal{W}^*_{ad}:=\Big\{ v\in L^\infty(\overline{Q_T})~~a.s.: l^{*}_2(x,y,t) \leq v(x,y,t)\leq l^{**}_2(x,y,t),~~a.e.~~a.s. \Big\},
\end{equation*}
\begin{equation*}
\mathcal{C}^{1*}_{ad}:=\Big\{ v\in L^\infty(\overline{Q_T})~~a.s.: l^{*}_3(x,y,t) \leq v(x,y,t)\leq l^{**}_3(x,y,t),~~a.e.~~a.s.\Big\},
\end{equation*}
\begin{equation*}
\mathcal{W}^{1*}_{ad}:=\Big\{  v\in L^\infty(\overline{Q_T})~~a.s.: l^{*}_4(x,y,t) \leq v(x,y,t)\leq l^{**}_4(x,y,t),~~a.e.~~a.s.\Big\},
\end{equation*}
 there exists a unique value ${J^*}(\beta_1,\beta_2,s_1, s_2)$ such that for every sequence $\Big\{(\beta^n_1,\beta^n_2,s^n_1, s^n_2)\Big\}_{n=0}^{\infty}$ in $\mathcal{C}_{ad}\times\mathcal{W}_{ad}\times\mathcal{C}^1_{ad}\times\mathcal{W}^1_{ad}$ converging to $(\beta_1,\beta_2,s_1, s_2)$ in 
 \[{L^2_4(\overline{Q_T})}:=L^2(\overline{Q_T})\times L^2(\overline{Q_T})\times L^2(\overline{Q_T})\times L^2(\overline{Q_T}),\]
  we have ${J}(\beta^n_1,\beta^n_2,s^n_1, s^n_2)\rightarrow{J^*}(\beta_1,\beta_2,s_1, s_2)$. Using the Lebesgue Dominated Convergence Theorem  \cite{26.}, we deduce that the function
\begin{equation}\label{AA*1ee}
\mathcal{J}(\beta_1,\beta_2,s_1, s_2)=
\begin{cases}
{J}(\beta_1,\beta_2,s_1, s_2),~~~(\beta_1,\beta_2,s_1, s_2)\in\mathcal{C}_{ad}\times\mathcal{W}_{ad}\times\mathcal{C}^1_{ad}\times\mathcal{W}^1_{ad},\\
{J^*}(\beta_1,\beta_2,s_1, s_2),~~(\beta_1,\beta_2,s_1, s_2)\in\mathcal{C}^*_{ad}\times\mathcal{W}^*_{ad}\times\mathcal{C}^{1*}_{ad}\times\mathcal{W}^{1*}_{ad}\setminus C_4(\overline{Q_T}),\\
+\infty,~~~~~~~~~~~~~~~~~(\beta_1,\beta_2,s_1, s_2)\in(\mathcal{C}^*_{ad}\times\mathcal{W}^*_{ad}\times\mathcal{C}^{1*}_{ad}\times\mathcal{W}^{1*}_{ad})^c\cap L^{1}_4({\overline{Q_T}}),\\
\end{cases}
\end{equation}
 is lower semicontinuous with respect to $(\beta_1,\beta_2,s_1, s_2)$ in ${L_4^{1}(\overline{Q_T})}$.
 Thus, from Theorem \ref{v principle}  (see Appendix),   we  conclude that for each positive $\epsilon$,  there exists $(\beta^\epsilon_1,\beta^\epsilon_2,s^\epsilon_1, s^\epsilon_2)$ such that
 \[\mathcal{J}(\beta^\epsilon_1,\beta^\epsilon_2,s^\epsilon_1, s^\epsilon_2)<\mathcal{J}(\beta_1,\beta_2,s_1, s_2)+\]
\begin{equation}\label{AA*1}
\epsilon^{\frac{1}{2}}\Big(\|\beta_1-\beta_1^\epsilon\|_{L^1(\overline{Q_T})}+\|\beta_2-\beta_2^\epsilon\|_{L^1(\overline{Q_T})}+\|s_1-s_1^\epsilon\|_{L^1(\overline{Q_T})}+\|s_2-s_2^\epsilon\|_{L^1(\overline{Q_T})}\Big),
\end{equation}
\[\forall(\beta_1,\beta_2,s_1, s_2)\not=(\beta^\epsilon_1,\beta^\epsilon_2,s^\epsilon_1, s^\epsilon_2).\]
In addition,  each sequence $\Big\{(\beta^n_1,\beta^n_2,s^n_1, s^n_2)\Big\}_{n=0}^{\infty}$ in $\mathcal{C}_{ad}\times\mathcal{W}_{ad}\times\mathcal{C}^1_{ad}\times\mathcal{W}^1_{ad}$ converging to $(\beta^\epsilon_1,\beta^\epsilon_2,s^\epsilon_1, s^\epsilon_2)$  in $L_4^2(\overline{Q_T})$ is a Cauchy sequence in $L_4^2(\overline{Q_T})$. If we assume that for each $n\in\mathbb{N}_0$, $(f_1^{n},f_2^{n})$ is the solution of   \eqref{YAS1111}--\eqref{YAS1*ekhco} corresponding to $(\beta^n_1,\beta^n_2,s^n_1, s^n_2)$, according to Theorem \ref{??}, the sequence $\Big\{(f_1^{n},f_2^{n})\Big\}_{n=0}^{\infty}$
  is a Cauchy sequence in $L^2_2(\overline{Q_T})$.
  So, using the Lebesgue Dominated Convergence Theorem  \cite{26.}, we derive that  for each $p>5$, there exist subsequences $\Big\{(\beta^{n_k}_1,\beta^{n_k}_2,s^{n_k}_1, s^{n_k}_2)\Big\}_{k=0}^{\infty}$ and $\Big\{(f_1^{n_k},f_2^{n_k})\Big\}_{k=0}^{\infty}$, which are Cauchy in 
\[L_4^p(\overline{Q_T}):= L^p(\overline{Q_T})\times L^p(\overline{Q_T})\times L^p(\overline{Q_T})\times L^p(\overline{Q_T}),\]
and
  \[L_2^p(\overline{Q_T}):= L^p(\overline{Q_T})\times L^p(\overline{Q_T}),\] 
respectively. Thus, using Lemma \ref{Lemma 1}  (see Appendix), it is easy to see that the sequences\newline
\[\Big\{(\mathfrak{f}_1^{n_k},\mathfrak{f}_2^{n_k})\Big\}_{k=0}^{\infty}=\Big\{(e^{-\sum_{i=0}^nh_{i,1}B_i}f_1^{n_k},e^{-\sum_{i=0}^nh_{i,2}B_i}f_2^{n_k})\Big\}_{k=0}^{\infty},\]
and 
\[\Big\{(\mathfrak{z}_{f_1^{n_k}},\mathfrak{z}_{f_2^{n_k}})\Big\}_{k=0}^{\infty}=\Big\{(e^{\sum_{i=0}^nh_{i,1}B_i}z_{f_1^{n_k}},e^{\sum_{i=0}^nh_{i,2}B_i}z_{f_2^{n_k}})\Big\}_{k=0}^{\infty},\]
 are Cauchy in
 \[W^{2,1,2}_p(\overline{Q_T}):= W^{2,1}_p(\overline{Q_T})\times W^{2,1}_p(\overline{Q_T}),\]
 where  $(z_{f_1^{n_k}},z_{f_2^{n_k}})$ is the solution of  \eqref{YAS*1e1}--\eqref{YAS1*e}  corresponding to  $(\beta^{n_k}_1,\beta^{n_k}_2,s^{n_k}_1, s^{n_k}_2)$.
Therefore,  we conclude that there exist \[A=(f_1^\epsilon,f_2^\epsilon,z^\epsilon_{f_1},z^\epsilon_{f_2},z^\epsilon_{g_1},z^\epsilon_{g_2},z^\epsilon_{g_3},z^\epsilon_{g_4}),\] and sequences
 \[\Big\{\smallunderbrace{(f_1^{m},f_2^{m},z_{f_1^{m}},z_{f_2^{m}},z_{g_1^{m}},z_{g_2^{m}},z_{g_3^{m}},z_{g_4^{m}})}_{A_m}\Big\}_{m=0}^{\infty},\]
 and
  \[\Big\{\smallunderbrace{(\beta^{m}_1,\beta^{m}_2,s^{m}_1, s^{m}_2)}_{B_m}\Big\}_{m=0}^{\infty},\]
 such that
\[\lim_{m\rightarrow\infty}A^*_m=A^*,~in~W^{2,1,4}_p(\overline{Q_T})\times L^p_4(\overline{Q_T}),~~\lim_{m\rightarrow\infty}B_m=(\beta^\epsilon_1,\beta^\epsilon_2,s_1^\epsilon,s_2^\epsilon),~in~L^p_4(\overline{Q_T}),\]
where
\[A_m^*=\]
\[(e^{-\sum_{i=0}^nh_{i,1}B_i}f_1^{m},e^{-\sum_{i=0}^nh_{i,2}B_i}f_2^{m},e^{\sum_{i=0}^nh_{i,1}B_i}z_{f_1^{m}},e^{\sum_{i=0}^nh_{i,2}B_i}z_{f_2^{m}},z_{g_1^{m}},z_{g_2^{m}},z_{g_3^{m}},z_{g_4^{m}}),\]
 and
  \[A^*=(e^{-\sum_{i=0}^nh_{i,1}B_i}f_1^\epsilon,e^{-\sum_{i=0}^nh_{i,2}B_i}f_2^\epsilon,e^{\sum_{i=0}^nh_{i,1}B_i}z^\epsilon_{f_1},e^{\sum_{i=0}^nh_{i,2}B_i}z^\epsilon_{f_2},z^\epsilon_{g_1},z^\epsilon_{g_2},z^\epsilon_{g_3},z^\epsilon_{g_4}).\]
Using the proof of Theorem \ref{ncon} and tangent-normal cone techniques \cite{25.}, for  $\epsilon$ small enough we have
\begin{equation}\label{YA1111}
|\beta_1^\epsilon-\mathcal{F}\Big(b_1-\dfrac{\partial z^\epsilon_{g_1}}{\partial t},l_1^{*},l_1^{**}\Big)|\leq \epsilon^{\frac{1}{2}},~|\beta_2^\epsilon-\mathcal{F}\Big(b_2-\dfrac{\partial z^\epsilon_{g_2}}{\partial t},l_2^{*},l_2^{**}\Big)|\leq\epsilon^{\frac{1}{2}},~a.e.~(x,y,t)\in\Omega\times [0,T],~a.s.
\end{equation}
\begin{equation}\label{YA1***}
|s_1^\epsilon-\mathcal{F}\Big(r_3^*-\dfrac{\partial z^\epsilon_{g_3}}{\partial t},l_3^*,l_3^{**}\Big)|\leq\epsilon^{\frac{1}{2}},
~|s_2^\epsilon-\mathcal{F}\Big(r_4^*-\dfrac{\partial z^\epsilon_{g_4}}{\partial t},l_4^*,l_4^{**}\Big)|\leq\epsilon^{\frac{1}{2}},~a.e.~(x,y,t)\in\Omega\times[0,T].~a.s.
\end{equation}
From  Theorem \ref{???}  (see Appendix) and the Gronwall inequality, it is clear that  the function
\[\mathcal{H}:\mathcal{C}_{ad}\times\mathcal{W}_{ad}\times\mathcal{C}^1_{ad}\times\mathcal{W}^1_{ad}\rightarrow\mathcal{C}_{ad}\times\mathcal{W}_{ad}\times\mathcal{C}^1_{ad}\times\mathcal{W}^1_{ad},\]
\[ \mathcal{H}(\beta_1,\beta_2,s_1, s_2)=\]
\[\left(\mathcal{F}\Big(b_1-\dfrac{\partial z_{g_1}}{\partial t},l_1^{*},l_1^{**}\Big),\mathcal{F}\Big(b_2-\dfrac{\partial z_{g_2}}{\partial t},l_2^{*},l_2^{**}\Big),\mathcal{F}\Big(r_3^*-\dfrac{\partial z_{g_3}}{\partial t},l_3^*,l_3^{**}\Big),\mathcal{F}\Big(r_4^*-\dfrac{\partial z_{g_4}}{\partial t},l_4^*,l_4^{**}\Big)\right),\]
where $z_{g_1}$, $z_{g_2}$, $z_{g_3}$ and $z_{g_4}$ are adjoint states corresponding to $(\beta_1,\beta_2,s_1, s_2)$, for almost every $\omega\in\bf{\Omega}$ has a unique fixed point $(\beta_1^*,\beta^*_2,s^*_1, s^*_2)$. Also using \eqref{AA*1} and Theorems  \ref{??} and \ref{???}, it is easy to show that  the sequence $\{(\beta_1^\epsilon,\beta^\epsilon_2,s^\epsilon_1, s^\epsilon_2)\}$ has a subsequence, which converges to $(\beta_1^*,\beta^*_2,s^*_1, s^*_2)$ in $L_4^2(\overline{Q_T})$, as $\epsilon$ converges to zero. Therefore, one can arrive at
\[\mathcal{J}(\beta_1^*,\beta^*_2,s^*_1, s^*_2)\leq\mathcal{J}(\beta_1,\beta_2,s_1, s_2),~~\forall(\beta_1,\beta_2,s_1, s_2)\in \mathcal{C}_{ad}\times\mathcal{W}_{ad}\times\mathcal{C}^1_{ad}\times\mathcal{W}^1_{ad},~~a.s.\]
which results in
\[{J}(\beta_1^*,\beta^*_2,s^*_1, s^*_2)\leq{J}(\beta_1,\beta_2,s_1, s_2),~~\forall(\beta_1,\beta_2,s_1, s_2)\in \mathcal{C}_{ad}\times\mathcal{W}_{ad}\times\mathcal{C}^1_{ad}\times\mathcal{W}^1_{ad}.~~a.s.\]
\qed

\section{Conclusions}\label{cremar}
 The model considered in this paper is a system of two stochastic nonlinear parabolic equations modelling the competition between two languages spoken in the same area. Since the evolution of the population of the speakers of each  language depends  on  several factors including the growth rate of the population and the status  of the language, the population with lower growth rate or status in compare with another one may be in danger of extinction. Clearly, by controlling the mentioned factors we can control the evolution of the populations. But, it is worth mentioning that changing these factors can be very costly and may have dangerous side effects on the population densities. These reasons encourage us to study the optimal control for the model of language competition to adopt more accurate and efficient  policies in order to control the evolution of the populations. For the reader's convenience, we briefly highlight our contributions as follows:
\begin{itemize}
\item We have added the stochastic terms to the model of competition between the languages to deal with a trustable model, which enables us to adopt more accurate policies. It is proved that the stochastic model has a unique solution.  
\item  We have considered four factors consisting of the growth rates of the populations and the status of the languages as the control variables. Then, by considering the cost function \eqref{focpcostf}, we have studied the optimal control problem for the model.
\item The stochastic adjoint equations \eqref{YAS*1e1}--\eqref{HM345**ee} are obtained, which help us to present the explicit forms of control variables and prove the existence and uniqueness of optimal control. It is shown the system of adjoint equations has a unique solution.   
\item By presenting some properties for the solutions of the problem \eqref{YAS1111}--\eqref{YAS1*ekhco} and the adjoint system \eqref{YAS*1e1}--\eqref{HM345**ee}   corresponding to different control variables in their admissible control sets (see Theorems \ref{??} and \ref{???} in Appendix) and employing  normal cone techniques and the adjoint equations, we have derived the necessary conditions for the stochastic optimal control variables. 
\item The Ekeland variational principle (Theorem \ref{v principle} in Appendix) is an effective tool in proving the existence
and uniqueness of stochastic optimal control variables.  This principle together with the properties of the solutions of the problem \eqref{YAS1111}--\eqref{YAS1*ekhco} and the adjoint system \eqref{YAS*1e1}--\eqref{HM345**ee}  (see Theorems \ref{??} and \ref{???} in Appendix) and some
useful results  enable us to prove that there exists unique optimal control.
\item  It is worthy of attention that we have presented the explicit forms of the stochastic optimal control
variables in Theorem \ref{ncon}.
\end{itemize}

\noindent{\textbf{Acknowledgments} }
\\
The authors are very grateful to reviewers for carefully reading the paper and for their valuable comments and 
suggestions which improved the original submission of this paper.

\begin{appendix}
 \section*{Appendix }
We provide here some essential mathematical concepts, lemmas and theorems which have been used in the
mathematical analysis throughout the paper.\\
\begin{defin}\label{121233}
\textup{\cite{steve345}} Let $B(t)$ be  1-dimensional Brownian motion  and $\mathfrak{F}(t)$ be an associated filtration. An Ito process is a stochastic process  of the form
\[X(t)=X(0)+\int_0^tu(s)ds+\int_0^tv(s)dB(s),\]
where $X(0)$ is nonrandom and $u(t)$ and $v(t)$ are adapted stochastic processes and
\[\mathbb{P}\Big[\int_0^t(v(s,\omega))^2ds<\infty ~\textup{for all}~ t\geq 0\Big]=1,~~\mathbb{P}\Big[\int_0^t|u(s,\omega)|ds<\infty ~\textup{for all}~ t\geq 0\Big]=1.\]

\end{defin}
\begin{theorem}
 \textup{\cite{steve345}} Let $f(t,x)$ be a function with continuous partial derivative $f_t$, $f_x$ and $f_{xx}$. Also assume that $X(t)$  be an Ito process as described in Definition \ref{121233}.
Then, for every $T>0$, we have
\[f(T,X(T))=f(0,X(0))+\int_0^T f_t(t,X(t))dt+\]
\[\int_0^T f_x(t,X(t))u(t)dt+\int_0^T f_x(t,X(t))v(t)dB(t)+\dfrac{1}{2}\int_0^T f_{xx}(t,X(t))v^2(t)dt.\]
\end{theorem}
\begin{theorem}
 \textup{\cite{oksendal}} Let $f(t,x)=(f_1(t,x),\cdots,f_n(t,x))$ be a $C^2$-smooth map from $[0,\infty)\times\mathbb{R}^n$ into $\mathbb{R}^m$ and 
$dX(t)=u(t)dt+v(t)dB(t)$ be an $n$-dimensional Ito process, which is of the following form
\begin{eqnarray} 
\left\lbrace \begin{array}{l} 
dX_1=u_1dt+v_{11}dB_{1}+\cdots+v_{1m}dB_{m},\\
\vdots\\
dX_n=u_ndt+v_{n1}dB_{1}+\cdots+v_{nm}dB_{m},
\end{array}\right. 
\end{eqnarray}
where $B=(B_1,\cdots,B_m)$ is an $m$-dimensional Brownian motion. 
Then, the process
\[Y(t)=f(t,X(t)),\]
is an Ito process and
\[dY_k(t)=\dfrac{\partial f_k}{\partial t}(t,X(t))dt+\sum_{i=0}^n\dfrac{\partial f_k}{\partial x_i}(t,X(t))dX_i(t)+\dfrac{1}{2}\sum_{i=0}^n\sum_{j=0}^n\dfrac{\partial^2 f_k}{\partial x_i\partial x_j}(t,X(t))dX_i(t)dX_j(t),\]
where $dB_idB_j=\delta_{ij}dt$ and $dtdt=dB_idt=dtdB_i=0$.
\end{theorem}
\begin{lem}
 \textup{\cite{steve345}} Let $X_1(t)$ and  $X_2(t)$ be Ito processes. Then
\[d(X_1(t)X_2(t))=X_1(t)dX_2(t)+X_2(t)dX_1(t)+dX_1(t)dX_2(t).\]
\end{lem}

\begin{theorem}\label{measure}
Let $(X,\Sigma)$ be a measurable space and $(\mathbb{R},\mathcal{B})$ be the Borel measurable space. Also, assume that $f:X\rightarrow\mathbb{R}$ and $g:X\rightarrow\mathbb{R}$ are measurable functions and $h:\mathbb{R}\longrightarrow\mathbb{R}$ is a continuous function, then
\begin{flushleft}
• $fg$ is measurable.
\end{flushleft}
\begin{flushleft}
• $h$ is measurable.
\end{flushleft}
\begin{flushleft}
• $hof$ is measurable if $f$ is finite.
\end{flushleft}
\begin{flushleft}
If $u=f$ $a.e.$ on $X$, then $u$ is a measurable function.
\end{flushleft}
\end{theorem}

\begin{defin} \textup{\cite{21.,23.} Let $0<\alpha<1$ and $\Omega\subset \mathbb{R}^n$ be bounded. Then $f\in C^{\alpha,\frac{\alpha}{2}}(\overline{\Omega}\times[0,T])$ if there exists a positive constant $C$ such that
\[|f(x_1,t_1)-f(x_2,t_2)|\leq C\Big(|x_1-x_2|^2+|t_1-t_2|\Big)^{\frac{\alpha}{2}},~\forall x_1,x_2\in\overline{\Omega},~\forall t_1,t_2\in[0,T].\] }
\textup{Furthermore, for any nonnegative integer $k$}
\[ C^{2k+\alpha,k+\frac{\alpha}{2}}(\overline{\Omega}\times[0,T]):=\{ f\in C^{\alpha,\frac{\alpha}{2}}(\overline{\Omega}\times[0,T]): \partial_x^{\beta} \partial_t^{i} f\in C^{\alpha,\frac{\alpha}{2}}(\overline{\Omega}\times[0,T]),~~   |\beta |+2i\leq 2k\}.\]
\end{defin}

\begin{theorem}\label{v principle}
\textup{\cite{24.}} Let $X$ be a complete metric space and let $f : X \rightarrow (-\infty, +\infty]$
be  lower semicontinuous and bounded from below and $\not\equiv+\infty$. Let $\epsilon > 0$ and $x_\epsilon\in X$ be such that
$f(x_\epsilon)\leq\inf\Big\{f(x): x\in X\Big\}+\epsilon.$
Then, there exists $y_\epsilon\in X$ such that
\[f(y_\epsilon)\leq f(x_\epsilon),~~~d(x_\epsilon,y_\epsilon)\leq\epsilon^{\frac{1}{2}},\]
\[f(y_\epsilon) < f(x) + \epsilon^{\frac{1}{2}}d(y_\epsilon, x),~~ \forall x \neq y_\epsilon.\]
\end{theorem}
\begin{defin} 
\textup{Let $Q_T=\Omega\times (0,T)$, then 
 we define 
\[W_{p}^{2,1 }(Q_T ):=\left\{u\in L^{p} (Q_T ) : \partial_x^\alpha \partial_t^k u\in L^{p} (Q_T ),~~   |\alpha |+2k\leq 2\right\},\]
with the norm
$||u||_{{W_{p}^{2,1}} (Q_T ) }:=\displaystyle{\sum _{|\alpha|+2k\leq 2}}||\partial_x^\alpha \partial
_t^k u||_{L^{p} (Q_T )}$.}
\end{defin}

\begin{lem}\label{Lemma 1} 
 Let $p>5$, $\Omega\subset \mathbb{R}^n$ be bounded with boundary $\partial\Omega\in C^2$ and
 \begin{equation}\label{re-diff1111}
\dfrac{\partial c}{\partial t}-\sum_{i=1}^n\sum_{j=1}^n a_{ij}(x,t)D_{ij}c+\sum_{i=1}^n  b_i(x,t)D_ic+d(x,t)c=g(x,t),~~(x,t)\in Q_T=\Omega\times (0,T),
\end{equation}
\[c(x,t)=\psi(x,t),~~(x,t)\in \partial Q_T\setminus(\Omega\times\{t=T\}),\]
where $\psi\in W^{2,1}_p(Q_T),~ a_{ij},b_i,d\in C(\overline{\Omega}\times[0,T])$, $a_{ij}=a_{ji}$ and for constants $0<C_1\leq C_2$
\[C_1|\lambda|^2\leq \sum_{i=1}^n\sum_{j=1}^n a_{ij}(x,t)\lambda_i\lambda_j\leq C_2|\lambda|^2,~~\forall\lambda\in \mathbb{R}^n,(x,t)\in Q_T.\]
Also assume that $g\in L^p(Q_T)$, then \eqref{re-diff1111} has a unique solution $c\in W^{2,1}_p(Q_T)$ such that
\begin{equation}
 {||c ||}_{W^{2,1}_p(Q_T)}\le  \mu\left({||c ||}_{L^p(Q_T)}+{||\psi ||}_{W^{2,1}_p(Q_T)}+||g ||_{L^p(Q_T)}\right),
\end{equation}
where $\mu$ depends on $T,~p, ~\Omega,~\|a_{ij}\|_{L^\infty(Q_T)},~\|b_{i}\|_{L^\infty(Q_T)}$ and $\|d\|_{L^\infty(Q_T)}$. Moreover, if $g\in C(\overline{Q_T})$, then
\begin{equation}
{||c ||}_{L^{\infty }(Q_T)}\le e^{\mu_0T}\Big(\|\psi\|_{L^{\infty }(\partial Q_T\setminus(\Omega\times\{t=T\})}+T||g ||_{L^\infty(Q_T)}\Big),
\end{equation} 
where $\mu_0=0$ if $d\geq 0$ and  $\mu_0=-\displaystyle \inf_{\overline{Q_{T}}}d$ otherwise.
\end{lem}
\noindent{\bf{Proof}} See \cite{21.,22.,Esmaili2}.
\qed

\begin{lem}\label{Lemma 1ito} Let $p>5$, $\Omega\subset \mathbb{R}^n$ be bounded with boundary $\partial\Omega\in C^2$ and
 \begin{equation}\label{ito3}
d c-\sum_{i=1}^n\sum_{j=1}^n a_{ij}(x,t)D_{ij}cdt+\sum_{i=1}^n  b_i(x,t)D_icdt+h(x,t)cdt=g(x,t)dt+c\sum_{k=0}^mh_k(x)dB_k,
\end{equation}
\[~~(x,t)\in Q_T=\Omega\times (0,T),\]
\[c(x,t)=\psi(x,t),~~(x,t)\in \partial Q_T\setminus(\Omega\times\{t=T\}),\]
\[\mathfrak{M}=\max\{\|a_{ij}\|_{L^\infty(Q_T)},~\|b_{i}\|_{L^\infty(Q_T)},~\|h\|_{L^\infty(Q_T)},~\|D_ih_k\|_{L^\infty(Q_T)},~\|D_{ij}h_k\|_{L^\infty(Q_T)}\},\]
where $\psi\in W^{2,1}_p(Q_T),~ a_{ij},b_i,h,D_{i,j}h_k\in C(\overline{\Omega}\times[0,T])$, $a_{ij}=a_{ji}$, and $h_k(x)=0$ on $\partial \Omega$, and for constants $0<C_1\leq C_2$
\begin{equation}\label{ghma}
C_1|\lambda|^2\leq \sum_{i=1}^n\sum_{j=1}^n a_{ij}(x,t)\lambda_i\lambda_j\leq C_2|\lambda|^2,~~\forall\lambda\in \mathbb{R}^n,(x,t)\in Q_T.
\end{equation}
Also assume that $g\in L^p(Q_T)$, then  for almost every $\omega\in \bf{\Omega}$, the problem \eqref{ito3} has a unique solution  $c(x,t )\in C(Q_T)$ such that $c(x,t )=e^{\sum_{k=0}^mh_k(x)B_k}C(x,t)$ and
\begin{equation}\label{mohem1}
 {||C||}_{W^{2,1}_p(Q_T)}\le  \mu(\omega)\left({||C ||}_{L^p(Q_T)}+{||\psi ||}_{W^{2,1}_p(Q_T)}+||e^{-\sum_{k=0}^mh_k(x)B_k}g ||_{L^p(Q_T)}\right),
\end{equation}
where
$\mu(\omega)$ depends on $T,~p, ~\Omega,$ and $\mathfrak{M}$. Moreover, if $g\in C(\overline{Q_T})$, then  $c$ is $\mathfrak{F}(t)-$adapted and for almost every $\omega\in\bf{\Omega}$, we have
\begin{equation}
{||c ||}_{L^{\infty }(Q_T)}\le {\mu'_0(\omega)}\Big(\|\psi\|_{L^{\infty }(\partial Q_T\setminus(\Omega\times\{t=T\})}+T||g ||_{L^\infty(Q_T)}\Big),
\end{equation} 
where  
 $\mu'_0(\omega)$ depends on $T,~p, ~\Omega,$ and $\mathfrak{M}$.  Also if $\psi=0$, then 
  \begin{equation}\label{mohem6}
  \int_{\Omega} (c(x,t))^2dx\leq e^{\lambda_2^*(\omega)t}\int_0^t\int_{\Omega}  (g(x,s))^2dxds,\end{equation}
where  
 $\lambda_2^*(\omega)$ depends on $T,~p, ~\Omega,$ and $\mathfrak{M}$.
\end{lem}
\noindent{\it Proof} Let $C$ be the solution of 
\[d C -\sum_{i=1}^n\sum_{j=1}^n a_{ij}(x,t)D_{ij}Cdt+\sum_{i=1}^n  \Big(\smallunderbrace{b_i(x,t)-2\sum_{j=1}^na_{ij}\sum_{k=0}^m B_kD_jh_k(x)}_{b_i'}\Big)D_iCdt+{\psi_1}Cdt=\]
\begin{equation}\label{ito2}
e^{-\sum_{k=0}^mh_k(x)B_k}g(x,t)dt, ~~(x,t)\in Q_T,
\end{equation} 
\[C(x,t)=\psi(x,t),~~(x,t)\in \partial Q_T\setminus(\Omega\times\{t=T\}),\]
where
\[\psi_1=h\left(x ,t \right)+\sum_{i=1}^nb_i(x,t)\sum_{k=0}^m B_kD_i {h_k(x)}+\dfrac{1}{2}\sum_{k=0}^m(h_k(x))^2-\]
\[e^{-\sum_{k=0}^mh_k(x)B_k} \sum_{i=1}^n\sum_{j=1}^n a_{ij}(x,t)D_{ij}(e^{\sum_{k=0}^mh_k(x)B_k}).\]
Now, we define $A_n$ and $\mathfrak{F}_n(t)$ as follows
\[A_n=\{\omega\in{\bf{\Omega}}: |B_k(s)|<n+1,~0\leq s\leq t,~~k=0,1,\cdots,m\}\cap\mathbb{A},\]
and
\[\mathfrak{F}_n(t)=\{A\cap A_n:~A\in\mathfrak{F}(t)\},\]
where
\[\mathbb{A}=\{\omega\in{\bf{\Omega}} :~B_K \textup{ is a countinuous function with respect to }t\} .\]
Clearly $B_k(t)\Big|_{ A_n}$ is a measurable function with respect to  $( A_n,\mathfrak{F}_n(t))$.
 Therefore, using Lemma \ref{Lemma 1} for almost every $\omega\in \bf{\Omega}$, the problem \eqref{ito2} has a unique solution such that $C(x,t )\in W^{2,1}_p(Q_T)$ and there exists  positive  $\mu_n$ 
 such that
\begin{equation}\label{ito4}
 {||C||}_{W^{2,1}_p(Q_T)}\le  \mu_n\left({||C ||}_{L^p(Q_T)}+{||\psi ||}_{W^{2,1}_p(Q_T)}+||e^{-\sum_{k=0}^mh_k(x)B_k}g ||_{L^p(Q_T)}\right),~~ \forall\omega\in A_n,
\end{equation}
where $\mu_n$ depends on $T,~p, ~\Omega,~\|a_{ij}\|_{L^\infty(Q_T)},~\|b'_{i}\|_{L^\infty(Q_T)}$ and $\|\psi_1\|_{L^\infty(Q_T)}$. Moreover, if $g\in C(\overline{Q_T})$, then there exists positive $\mu^n_0$ such that
\begin{equation}\label{ito6}
{||C ||}_{L^{\infty }(Q_T)}\le e^{\mu^n_0T}\Big(\|\psi\|_{L^{\infty }(\partial Q_T\setminus(\Omega\times\{t=T\})}+T||g ||_{L^\infty(Q_T)}\Big),~~\forall\omega\in A_n.
\end{equation} 
It is also clear that the solution of \eqref{ito2}, $C$, is a function of $B_k$ $(k=0,\cdots,m)$. From \eqref{ito2} and \eqref{ito6},  $C\Big|_{A_n}=\mathfrak{G}(B_0\Big|_{A_n},\cdots,B_m\Big|_{A_n})$ is a continuous function with respect to $B_k\Big|_{A_n}$ $(k=0,\cdots,m)$. Therefore, from Theorem \ref{measure},  $C\Big|_{A_n}$ is $\mathfrak{F}_n(t)-$adapted. Thus $C\Big|_{A_n}$ is $\mathfrak{F}(t)-$adapted. Also,   $C$ the solution of \eqref{ito2} can be defined (almost surely) as follows
\[C=C\Big|_{A_n},~~ \omega\in A_n, ~~n\in\mathbb{N}_0.\]
Therefore, for every borel set $B$, $C^{-1}(B)=\displaystyle\cup_{n=0}^\infty C\Big|_{A_n}^{-1}(B)\in\mathfrak{F}(t).$ So, $C$ is  $\mathfrak{F}(t)-$adapted.  \\
 On the other hand, we have
\[d(C^2)=2CdC=2C\sum_{i=1}^n\sum_{j=1}^n a_{ij}(x,t)D_{ij}Cdt-2C\sum_{i=1}^n  {b_i'}D_iCdt-2{\psi_1}C^2dt+2e^{-\sum_{k=0}^mh_k(x)B_k}g(x,t)Cdt,\]
\[ ~~(x,t)\in Q_T,\]
\[C(x,t)=\psi(x,t),~~(x,t)\in \partial Q_T\setminus(\Omega\times\{t=T\}).\]
Therefore if $\psi=0$, then
\[\int_{\Omega} d(C^2)dx=\]
\[\int_{\Omega}\Big(2C\sum_{i=1}^n\sum_{j=1}^n a_{ij}(x,t)D_{ij}C-2C\sum_{i=1}^n  {b_i'(x,t)}D_iC-2{\psi_1(x,t)}C^2+2e^{-\sum_{k=0}^mh_k(x)B_k}g(x,t)C\Big)dtdx=\]
\[\int_{\Omega}\Big(-2\sum_{i=1}^n\sum_{j=1}^n a_{ij}(x,t)D_{i}CD_jC-2C\sum_{i=1}^n  {b_i'(x,t)}D_iC-2{\psi_1(x,t)}C^2+2e^{-\sum_{k=0}^mh_k(x)B_k}g(x,t)C\Big)dtdx.\]
Hence, from \eqref{ghma}, we conclude that if $\psi=0$, then
\[\int_{\Omega} (C(x,t))^2dx=\]
\[\int_0^t\int_{\Omega}\Big(-2\sum_{i=1}^n\sum_{j=1}^n a_{ij}(x,s)D_{i}CD_jC-2C\sum_{i=1}^n  {b_i'(x,s)}D_iC-2{\psi_1(x,s)}(C)^2+2e^{-\sum_{k=0}^mh_k(x)B_k}g(x,s)C\Big)dsdx\leq\]
\[\int_0^t\int_{\Omega}\lambda^*(\omega)C^2dxds+\int_0^t\int_{\Omega}  e^{-2\sum_{k=0}^mh_k(x)B_k}(g(x,s))^2dxds,~~a.s.\]
where  
 $\lambda^*(\omega)$ depends on $T,~p, ~\Omega$ and $\mathfrak{M}$. Employing Gronwall inequality, we deduce that 
  \begin{equation}\label{eqG}
  \int_{\Omega} (C(x,t))^2dx\leq e^{\lambda_1^*(\omega)t}\int_0^t \int_{\Omega}  (g(x,s))^2dxds,~~a.s.\end{equation}
where  
 $\lambda_1^*(\omega)$ depends on $T,~p, ~\Omega,$ and $\mathfrak{M}$. \\
Since $C$ is $\mathfrak{F}(t)-$adapted, $c$ is $\mathfrak{F}(t)-$adapted. Thus, using the Ito product rule for $c=e^{\sum_{k=0}^mh_k(x)B_k}C$, we have
\[d(\smallunderbrace{e^{\sum_{k=0}^mh_k(x)B_k}C}_{c})=d(e^{\sum_{k=0}^mh_k(x)B_k})C+e^{\sum_{k=0}^mh_k(x)B_k}dC=\]
\begin{equation}\label{ito1}
Ce^{\sum_{k=0}^mh_k(x)B_k}\sum_{k=0}^mh_k(x)dB_k+\dfrac{{C}}{2}e^{\sum_{k=0}^mh_k(x)B_k}\sum_{k=0}^m(h_k(x))^2dt+e^{\sum_{k=0}^mh_k(x)B_k}d{C}.\end{equation}
Thus, from \eqref{ito2} and \eqref{ito1}, we can conclude that
\begin{equation}\label{ito7}
d c-\sum_{i=1}^n\sum_{j=1}^n a_{ij}(x,t)D_{ij}cdt+\sum_{i=1}^n  b_i(x,t)D_icdt+h(x,t)cdt=g(x,t)dt+c\sum_{k=0}^mh_kdB_k,~~(x,t)\in Q_T,
\end{equation}
\[c(x,t)=\psi(x,t),~~(x,t)\in \partial Q_T\setminus(\Omega\times\{t=T\}).\]
Hence, from \eqref{ito6} and \eqref{ito7}, we conclude that for almost every $\omega\in\bf{\Omega}$, \eqref{ito3} has a continuous solution $c$. Moreover, if $g\in C(\overline{Q_T})$, then 
\begin{equation}
{||c ||}_{L^{\infty }(Q_T)}\le \mu'_0(\omega)\Big(\|\psi\|_{L^{\infty }(\partial Q_T\setminus(\Omega\times\{t=T\})}+T||g ||_{L^\infty(Q_T)}\Big),~~a.s.
\end{equation} 
where  
 $\mu'_0(\omega)$ depends on $T,~p, ~\Omega,$ and $\mathfrak{M}$. Also if $\psi=0$, then from \eqref{eqG} we have
   \begin{equation}
  \int_{\Omega} (c(x,t))^2dx\leq e^{\lambda_2^*(\omega)t}\int_0^t\int_{\Omega}  (g(x,s))^2dxds,~~a.s.\end{equation}
where  
 $\lambda_2^*(\omega)$ depends on $T,~p, ~\Omega$ and $\mathfrak{M}$.
  Now, we want to show that the problem \eqref{ito3} has a unique solution. Let ${C}^1$ and ${C}^2$ be the solutions of \eqref{ito3}. Thus, $U=C^1-C^2$ is the solution of
\begin{equation}
d U-\sum_{i=1}^n\sum_{j=1}^n a_{ij}(x,t)D_{ij}Udt+\sum_{i=1}^n  b_i(x,t)D_iUdt+h(x,t)Udt=U\sum_{k=0}^mh_kdB_k,~~(x,t)\in Q_T,
\end{equation}
\[U(x,t)=0,~~(x,t)\in \partial Q_T\setminus(\Omega\times\{t=T\}).\]
On the other hand, using the Ito product rule for $U_1=e^{-\sum_{k=0}^mh_k(x)B_k}U$, one can conclude that 
\[dU_1=e^{-\sum_{k=0}^mh_k(x)B_k}dU+Ude^{-\sum_{k=0}^mh_k(x)B_k}+de^{-\sum_{k=0}^mh_k(x)B_k}dU=\]
\[e^{-\sum_{k=0}^mh_k(x)B_k}dU-Ue^{-\sum_{k=0}^mh_k(x)B_k}\sum_{k=0}^mh_k(x)dB_k+\dfrac{{U}}{2}e^{-\sum_{k=0}^mh_k(x)B_k}\sum_{k=0}^m(h_k(x))^2dt-\]
\[Ue^{-\sum_{k=0}^mh_k(x)B_k}\sum_{k=0}^m(h_k(x))^2dt=\]
\[e^{-\sum_{k=0}^mh_k(x)B_k}dU-U_1\sum_{k=0}^mh_k(x)dB_k-\dfrac{{U_1}}{2}\sum_{k=0}^m(h_k(x))^2dt.\]
Therefore, $U_1=e^{-\sum_{k=0}^mh_k(x)B_k}U$ is the solution of
\[d U_1 -\sum_{i=1}^n\sum_{j=1}^n a_{ij}(x,t)D_{ij}U_1dt+\sum_{i=1}^n  \Big(b_i(x,t)-2\sum_{j=1}^na_{ij}\sum_{k=0}^m B_kD_jh_k(x)\Big)D_iU_1dt+\]
\[\Big(h\left(x ,t \right)+\sum_{i=1}^nb_i(x,t)\sum_{k=0}^m B_kD_i {h_k(x)}+\dfrac{1}{2}\sum_{k=0}^m(h_k(x))^2-\]
\[e^{-\sum_{k=0}^mh_k(x)B_k} \sum_{i=1}^n\sum_{j=1}^n a_{ij}(x,t)D_{ij}(e^{\sum_{k=0}^mh_k(x)B_k} )\Big)U_1dt=0,\] 
\[U_1(x,t)=0,~~(x,t)\in \partial Q_T\setminus(\Omega\times\{t=T\}).\]
Finally, using Lemma \ref{Lemma 1}, one can deduce that   $U_1=0$ for almost every $\omega\in \bf{\Omega}$, which results in the uniqueness of the solution of the problem \eqref{ito3} for almost every $\omega\in\bf{\Omega}$.
\qed

\begin{theorem}\label{??}
Let $p>5$ and  $({f_1}^1,{f_2}^1)$ and $({f_1}^2,{f_2}^2)$  be the exact solutions of  the problem \eqref{YAS1111}--\eqref{YAS1*ekhco} corresponding to $(\beta^1_1,s^1_1,\beta^1_2,s^1_2)\in \mathcal{C}_{ad}\times\mathcal{C}^1_{ad}\times\mathcal{W}_{ad}\times\mathcal{W}^1_{ad}$ and $(\beta^2_1,s^2_1,\beta^2_2,s^2_2)\in \mathcal{C}_{ad}\times\mathcal{C}^1_{ad}\times\mathcal{W}_{ad}\times\mathcal{W}^1_{ad}$,  respectively.
 Then, for almost every $\omega\in\mathbf{\Omega}$, there exist positive   $\iota_0(\omega)$ and $\iota(\omega)$ such that 
\[{\|({f_1}^1,{f_2}^1)-({f_1}^2,{f_2}^2)\|}^p_{L^{\infty}(\Omega) }\le\]
\[ \iota_0(\omega)\int_0^t{\|( \beta^1_1,\beta^1_2,s^1_1, s^1_2)-(\beta^2_1,\beta^2_2,s^2_1, s^2_2) \|}^p_{L^{\infty}(\Omega) }ds,\]
and
\[{\|({f_1}^1,{f_2}^1)-({f_1}^2,{f_2}^2)\|}^2_{L^2(Q_t)}\leq\]
\[\iota(\omega)\int_0^t{\|(\beta^1_1,\beta^1_2,s^1_1, s^1_2)-(\beta^2_1,\beta^2_2,s^2_1, s^2_2) \|}^2_{L^2(Q_s)}ds,~~Q_t=\Omega\times(0,t),\]
where $\iota_0(\omega)$ and $\iota(\omega)$ depend on $~p, ~\Omega,~\|D_ih_{k,l}\|_{L^\infty(Q_T)},~\|D_{ij}h_{k,l}\|_{L^\infty(Q_T)}$ and $\|h_{k,l}\|_{L^\infty(Q_T)}$ $(k=0,1,~\cdots,n,l=1,2)$.
\end{theorem}
\noindent{\bf{ Proof}}  From the problem \eqref{YAS1111}--\eqref{YAS1*ekhco}, for $\dot{f}_1=f_1^1-f_1^2$, we have
\[d{\dot{f}_1}=\nabla.(D_{1}\nabla{  {\dot{f}_1}})dt-\nabla.(F_1\dot{f}_1)dt+\]
\[\Big(\smallunderbrace{R_1(f^1_1,f^1_2,\beta^1_1)+R(f^1_1,f^1_2,s^1_1,s^1_2)-R_1(f^2_1,f^2_2,\beta^2_1)-R(f^2_1,f^2_2,s^2_1,s^2_2)}_{g_1^*}\Big)dt+\]
\begin{equation}\label{mohem7}
\dot{f}_1\sum_{k=0}^n h_{k,1}(x,y)dB_k,~on~\Omega,\ t {\rm >}0,
\end{equation}
\[{\dot{f}_1}\left({x,y,}t \right){\rm =}0, ~on~\partial\Omega,~ {t >}0,~{\dot{f}_1}\left(x,y ,0\right){\rm =}0,~on~\Omega.\]
Using \eqref{mohem1}, one can deduce that 
$\dot{f}_1(x,y,t )=e^{\sum_{k=0}^nh_{k,1}(x,y)B_k}\dot{F}$ and
\begin{equation}
 {||\dot{F}||}_{W^{2,1}_p(Q_t)}\le  \mu^*(\omega)\left({||\dot{F} ||}_{L^p(Q_t)}+||g_1^* ||_{L^p(Q_t)}\right),
\end{equation}
where
$\mu^*(\omega)$ depends on $\|D_ih_{k,1}\|_{L^\infty(Q_T)},~\|D_{ij}h_{k,1}\|_{L^\infty(Q_T)}$ and $\|h_{k,1}\|_{L^\infty(Q_T)}$. Using t-Anisotropic Embedding Theorem \cite{21.}, one can conclude that

\begin{equation}\label{mohem2}
 (\dot{f}_1(x,y,t))^p\le  \mu_1^*(\omega)\left({||\dot{f}_1(x,y,s) ||}^p_{L^p(Q_t)}+||g_1^*(x,y,s) ||^p_{L^p(Q_t)}\right).~~a.s.
\end{equation} 
Similarly, it can be shown that
\begin{equation}\label{mohem3}
 (\dot{f}_2(x,y,t))^p\le  \mu_2^*(\omega)\left({||\dot{f}_2(x,y,s) ||}^p_{L^p(Q_t)}+||g_2^*(x,y,s) ||^p_{L^p(Q_t)}\right),~~a.s.
\end{equation} 
where
\[\dot{f}_2=f_2^1-f_2^2,\]
and
\[g_2^*=R_2(f^1_1,f^1_2,\beta^1_2)-R(f^1_1,f^1_2,s^1_1,s^1_2)-R_2(f^2_1,f^2_2,\beta^2_2)+R(f^2_1,f^2_2,s^2_1,s^2_2).\]
From \eqref{mohem2} and \eqref{mohem3}, we have
\[ \|\dot{f}_1(x,y,t)\|^p_{L^\infty(\Omega)}+\|\dot{f}_2(x,y,t)\|^p_{L^\infty(\Omega)}\le \]
\[ \mu_3^*(\omega)\Big({\int_0^t\|\dot{f}_1(x,y,s)\|^p_{L^\infty(\Omega)}+\|\dot{f}_2(x,y,s)\|^p_{L^\infty(\Omega)}ds}+\]
\[\int_0^t\|( \beta^1_1,\beta^1_2,s^1_1, s^1_2)-(\beta^2_1,\beta^2_2,s^2_1, s^2_2) \|^p_{L^{\infty }(\Omega)}ds\Big).~a.s.\]
Thus, applying the Gronwall inequality results in
\[ \|\dot{f}_1(x,y,t)\|_{L^\infty(\Omega)}+\|\dot{f}_2(x,y,t)\|_{L^\infty(\Omega)}\le \]
\[ \iota^*_0(\omega)\left(\int_0^t\|( \beta^1_1,\beta^1_2,s^1_1, s^1_2)-(\beta^2_1,\beta^2_2,s^2_1, s^2_2) \|^p_{L^{\infty }(\Omega)}ds\right)^{\frac{1}{p}}.~~a.s.\]
Also, from \eqref{mohem6}, \eqref{mohem7} and the Gronwall inequality, one can deduce that
\[{\|({f_1}^1,{f_2}^1)-({f_1}^2,{f_2}^2)\|}^2_{L^2(Q_t)}\leq\iota(\omega)\int_0^t{\|(\beta^1_1,\beta^1_2,s^1_1, s^1_2)-(\beta^2_1,\beta^2_2,s^2_1, s^2_2) \|}^2_{L^2(Q_s)}ds.~~a.s.\]
\qed

\begin{theorem}\label{???}
Let $p>5$ and   $(z_{f_{1}}^1,z_{f_{2}}^1,z_{g_{1}}^1,z_{g_{2}}^1,z_{g_{3}}^1,z_{g_{4}}^1)$ and $(z_{f_{1}}^2,z_{f_{2}}^2,z_{g_{1}}^2,z_{g_{2}}^2,z_{g_{3}}^2,z_{g_{4}}^2)$ be the exact solutions of the stochastic adjoint system \eqref{YAS*1e1}--\eqref{HM345**ee}  corresponding to $(\beta^1_1,s^1_1,\beta^1_2,s^1_2)\in \mathcal{C}_{ad}\times\mathcal{C}^1_{ad}\times\mathcal{W}_{ad}\times\mathcal{W}^1_{ad}$ and $(\beta^2_1,s^2_1,\beta^2_2,s^2_2)\in \mathcal{C}_{ad}\times\mathcal{C}^1_{ad}\times\mathcal{W}_{ad}\times\mathcal{W}^1_{ad}$,  respectively. Then, for almost every $\omega\in\mathbf{\Omega}$, there exists positive   $\kappa^*(\omega)$ such that 
\[{\|(z_{f_{1}}^1,z_{f_{2}}^1,z_{g_{1}}^1,z_{g_{2}}^1,z_{g_{3}}^1,z_{g_{4}}^1)-(z_{f_{1}}^2,z_{f_{2}}^2,z_{g_{1}}^2,z_{g_{2}}^2,z_{g_{3}}^2,z_{g_{4}}^2)\|}^p_{L^{\infty}(\Omega) }\le\]
\[ \kappa^*(\omega)\Big(\int_t^T\int_0^s{\|( \beta^1_1,\beta^1_2,s^1_1, s^1_2)-(\beta^2_1,\beta^2_2,s^2_1, s^2_2) \|}^p_{L^{\infty}(\Omega) }dlds+\]
\[\int_t^T{\|( \beta^1_1,\beta^1_2,s^1_1, s^1_2)-(\beta^2_1,\beta^2_2,s^2_1, s^2_2) \|}^p_{L^{\infty}(\Omega) }ds\Big),\]
where $\kappa^*(\omega)$  depends on $~p, ~\Omega,~\|D_ih_{k,l}\|_{L^\infty(Q_T)},~\|D_{ij}h_{k,l}\|_{L^\infty(Q_T)}$ and $\|h_{k,l}\|_{L^\infty(Q_T)}$ $(k=0,1,~\cdots,n,l=1,2)$.
\end{theorem}
\noindent{\bf{Proof}} Similar to the proof of Theorem \ref{??}, we can get the result.
\qed

\end{appendix}

\end{document}